\documentclass[preprint,authoryear, 12pt]{elsarticle}
\usepackage{setspace}
\usepackage{graphics}
\usepackage{subfigure}
\usepackage{psfrag}
\usepackage{epsfig}
\usepackage{pstricks}
\usepackage{amsthm}
\usepackage{amsmath}
\usepackage{amsfonts}
\usepackage{amssymb}
\usepackage{ae}
\usepackage{graphicx}
\usepackage{natbib}

\newtheorem{theorem}{Theorem}[section]

\theoremstyle{definition}

\newenvironment{frcseries}{\fontfamily{frc}\selectfont}{}
\newcommand{\textfrc}[1]{{\frcseries#1}}

\theoremstyle{remark}
\newtheorem{remark}[theorem]{Remark}

\numberwithin{equation}{section}

\biboptions{round}

\begin{document}
\begin{frontmatter}
\title{Nonlinear Neutral Inclusions: Assemblages of Coated Ellipsoids.}

\author{Silvia Jim\'enez\corref{cor1}}
\ead{sjimenez@colgate.edu}

\author{Bogdan Vernescu}
\ead{vernescu@wpi.edu}

\address{*Dept. of Mathematics, Colgate University, 13 Oaks Drive, Hamilton, NY 13346\corref{cor1} \\Dept. of Mathematical Sciences, Worcester Polytechnic Institute, 100 Institute Road, Worcester, MA 01609}

\cortext[cor1]{Corresponding author}


\begin{abstract}
The problem of determining nonlinear neutral inclusions in (electrical or thermal) conductivity is considered. Neutral inclusions, inserted in a matrix containing a uniform applied electric field, do not disturb the field outside the inclusions. The well known Hashin coated sphere construction is an example of a neutral inclusion. In this paper, we consider the problem of constructing neutral inclusions from nonlinear materials. In particular, we discuss assemblages of coated ellipsoids.  
\end{abstract}
\end{frontmatter}
\section{Introduction}
A neutral inclusion, when inserted in a matrix containing a uniform applied electric field, does not disturb the outside field.  Mansfield was the first to observe that  reinforced holes, ``neutral holes'', could be cut out of a uniformly stressed plate without disturbing the surrounding stress field in the plate \cite{Mansfield1953}.  
	
The well known Hashin coated sphere construction \cite{Hashin1962} is an example of a neutral coated inclusion for the conductivity problem. In \cite{Hashin1962b,Hashin1963} an exact expression for the effective conductivity of the coated sphere assemblage was found, which coincides with the Maxwell \cite{Maxwell1873} approximate formula. Thus the approximate formula is realizable and was shown to be an attainable bound for the effective conductivity of a composite, given the volume fractions of the two materials.  This construction was extended to coated confocal ellipsoids in \cite{Milton1981}.  Ellipsoids as neutral inclusions have been also studied in \cite{Kerker1975}.  Spheres and ellipsoids are not the only possible shapes for neutral inclusions; indeed in \cite{Milton2001} other shapes of neutral inclusions are constructed.

The existence of neutral inclusions was also found in the case of materials with imperfect interfaces, for which the potential (or displacement) field has discontinuities across these interfaces. For these materials neutral  inclusions have been studied in \cite{Lipton1996}, \cite{Benveniste1999} for the conductivity problem, in \cite{Lipton1997}, \cite{Lipton1997-2} for highly conducting interfaces, in \cite{Lipton1995}, \cite{Lipton1996-2}, \cite{Ru1998} for the elasticity problem, and for nonlinear materials in \cite{Lipton1999}.   

For other references related to neutral inclusions in composites see also \cite{Milton2002} and \cite{Vernescu2010} and the references therein.  

We consider here nonlinear materials for which the constitutive law relating the current $J$ to the electric field $\nabla u$  is described by a nonlinear constitutive model of the form
$$J=\sigma_1\left|\nabla u\right|^{p-2}\nabla u,$$
here $u$ is the potential, and $\sigma_1\left|\nabla u\right|^{p-2}$ is a nonlinear conductivity.
This constitutive model is used to describe the nonlinear behavior of several materials including nonlinear dielectrics \cite{Bueno2008,Garroni2001,Garroni2003,Kohn1998,Talbot1994,Talbot1994-2}, and is also used to model thermo-rheological and electro-rheological fluids \cite{Ruzicka2000,Antontsev2006,Ruzicka2008}, viscous flows in glaciology \cite{Glowinski2003}, and also in plasticity problems \cite{Atkinson1984,Suquet1993,PonteCastaneda1997,PonteCastaneda1999,Idiart2008}.  

In this paper we show that even for nonlinear materials, one can construct neutral inclusions by a suitable coating with a linear material.  In particular, we show that that a coated ellipsoid with core of phase $1$ (nonlinear material) surrounded by a coating of phase $2$ (linear material) can be constructed as a neutral inclusion.  In \cite{Jimenez2013}, we showed that coated spheres with nonlinear core and linear coating can be constructed as neutral inclusions.  

Since the equations for conductivity are local equations, one could continue to add similar aligned coated ellipsoids of various sizes without disturbing the prescribed uniform applied field surrounding the inclusions.  In fact, one can fill the entire space (aside from a set of measure zero) with assemblages of these aligned coated ellipsoids by adding coated ellipsoids of various sizes ranging to the infinitesimal and it is assumed that they do not overlap the boundary of the unit cell of periodicity.  The ellipsoids can be of any size, but the volume fraction $\theta_1$ (\ref{volumefraction}) of nonlinear material is the same for all ellipsoids.  While adding the coated ellipsoids, the flux of current and electrical potential at the boundary of the unit cell remains unaltered.  Therefore, the effective conductivity does not change.  

This paper is structured as follows: Section~\ref{Stat-Ellipsoids} provides the statement of the problem and the main result for an assemblage of coated ellipsoids and Section~\ref{Results} provides the proofs of the statements in Section~\ref{Stat-Ellipsoids}.  

\section{Assemblage of Coated Ellipsoids: Statement of the Problem.}  
\label{Stat-Ellipsoids}

We need to introduce ellipsoidal coordinates $\rho$, $\mu$, and $\nu$, which are defined implicitly as the solution of the set of equations \cite{Landau1984,Kellogg1953}

\begin{equation*}
\begin{cases}
\displaystyle \frac{x_{1}^2}{c_{1}^2+\rho}+\frac{x_{2}^2}{c_{2}^2+\rho}+\frac{x_{3}^2}{c_{3}^2+\rho}=1: 
\text{Confocal Ellipsoids}\\
\displaystyle \frac{x_{1}^2}{c_{1}^2+\mu}+\frac{x_{2}^2}{c_{2}^2+\mu}+\frac{x_{3}^2}{c_{3}^2+\mu}=1: \text{Hyperboloids of one sheet}\\
\displaystyle\frac{x_{1}^2}{c_{1}^2+\nu}+\frac{x_{2}^2}{c_{2}^2+\nu}+\frac{x_{3}^2}{c_{3}^2+\nu}=1: \text{Hyperboloids of two sheets}\\
\end{cases}
\end{equation*}

subject to the restrictions $$\rho>-c_{1}^2>\mu>-c_{2}^2>\nu>-c_{3}^2,$$ where $c_{1}$, $c_{2}$, and $c_{3}$ are fixed positive constants that determine the coordinate system, all confocal with the ellipsoid $$\displaystyle \frac{x_{1}^2}{c_{1}^2}+\frac{x_{2}^2}{c_{2}^2}+\frac{x_{3}^2}{c_{3}^2}=1.$$  One surface of each of the three families passes through each point in space, and the three surfaces are orthogonal.  The equations can be solved explicitly for the Cartesian coordinates in terms of the ellipsoidal coordinates.  For all permutations $j$, $k$, $l$ of $1$, $2$, $3$ we have
\begin{equation}
	\displaystyle x_{j}^2=\frac{(c_{j}^2+\rho)(c_{j}^2+\mu)(c_{j}^2+\nu)}{(c_{j}^2-c_{k}^2)(c_{j}^2-c_{l}^2)}.
\end{equation}  

The coordinate $\rho$ plays the role that the radius plays in spherical coordinates.  Our prototype ellipsoid is defined by the region $\rho<\rho_{e}$ with a nonlinear core $0<\rho<\rho_c$ and a linear coating $\rho_{c}<\rho<\rho_{e}$.  Within the ellipsoid the conductivity depends only on the coordinate $\rho$. 

We introduce the lengths 
\begin{equation}
\label{semiaxis}
	\displaystyle \textfrc{l}_{cj}=\sqrt{c_{j}^2+\rho_c} \text{, \hspace{3mm}$\textfrc{l}_{ej}=\sqrt{c_{j}^2+\rho_e}$, \hspace{3mm}$j=1,2,3$,}
\end{equation}
which represent the semi-axis lengths of the core and exterior surfaces of the coated ellipsoid, the volume fraction
\begin{equation}
		\label{volumefraction}
\theta_{1}=\frac{\textfrc{l}_{c1}\textfrc{l}_{c2}\textfrc{l}_{c3}}{\textfrc{l}_{e1}\textfrc{l}_{e2}\textfrc{l}_{e3}},
	\end{equation}
occupied by phase $1$ (nonlinear material in the core) and $\theta_{2}=1-\theta_{1}$, the volume fraction occupied by phase $2$ (linear material in the coating).

The coated ellipsoid is embedded in a medium with isotropic conductivity tensor $\sigma_1^*\textbf{I}$, where the value of $\sigma_1^*$ needs to be chosen so that the conductivity equations have a solution with the uniform field aligned in the $x_1$ direction in the region exterior to the ellipsoid.  Once this is done, it follows by the usual argument that $\sigma_1^*$ represents the effective conductivity in the $x_1$ direction of the assemblage of aligned ellipsoids, each identical within a scale factor to the given prototype.  We apply a linear electric field $\mathbf{E}\cdot\mathbf{x}=Ex_{1}$ at infinity, (where for simplicity $\mathbf{E}=E\mathbf{e^{1}}$, with $\mathbf{e^{1}}=(1,0,0)$ and $\mathbf{x}=(x_{1},x_{2},x_{3})$). 

Thus the problem of finding a neutral inclusion reduces to finding the electric potential $u$ that solves 
	\begin{equation}
		\label{PDENeutralInclusions}
		\begin{cases}
			\nabla\cdot\left(\sigma_{1}\left|\nabla u\right|^{p-2}\nabla u\right)=0 \;\;\text{in the core,}\\
			\nabla\cdot\left(\sigma_{2}\nabla u\right)=0 \;\;\text{in the coating,}
		\end{cases}
	\end{equation}
where the material conductivities are $\sigma_1\left|\nabla u\right|^{p-2}$ in the core, and $\sigma_2$ in the coating, with $\infty>\sigma_{1}>\sigma_{2}>0$, and satisfies continuity conditions  of the electric potential and of the normal component of the current at the interfaces.  

\section{Assemblage of Coated Ellipsoids: Results}
\label{Results}
 
Inside the coated ellipsoid, we ask that 
\begin{equation}
	\label{PDENeutralInclusions2}
	\begin{cases}
		\sigma_{1}\Delta_{p}u=0 & \mbox{ for } 0<\rho<\rho_{c} \\
		\sigma_{2}\Delta u=0 & \mbox{ for } \rho_c<\rho<\rho_e, \\
	\end{cases}
\end{equation}
where $\Delta_{p}u=\nabla\cdot(\left|\nabla u\right|^{p-2}\nabla u)$ represents the $p$-Laplacian ($p>1$), $\sigma_{1}$ and $\sigma_{2}$ are positive, together with the usual continuity conditions  of the electric potential and of the normal component of the current across the interfaces: 
\begin{equation}
	\label{conta}
	\text{$u$ continuous across $\rho=\rho_c$,}
\end{equation}
\begin{equation}
	\label{extbound}
	\text{$u=Ex_{1}$ at $\rho=\rho_e$,}
\end{equation}
and
\begin{equation}
	\label{transa}
	\text{$\sigma_{1}\left|\nabla u\right|^{p-2}\nabla u \cdot\mathbf{n}=\sigma_{2}\nabla u\cdot\mathbf{n}$,  across $\rho=\rho_c$,}
\end{equation}
\begin{equation}
	\label{transb}
	\text{$\sigma_{2}\nabla u\cdot \mathbf{n}=\sigma_{1}^*\nabla u \cdot\mathbf{n}$,  across $\rho=\rho_e$.}
\end{equation}

We look for a solution $u$ of (\ref{PDENeutralInclusions2}) of the form 
\begin{equation}
	\label{u-ni}
	u=\begin{cases}
		A_1x_{1} \text{     for $0<\rho<\rho_{c}$},\\
		\varphi(\rho)x_{1} \text{     for $\rho_{c}\leq\rho\leq\rho_{e}$}.
	\end{cases}
\end{equation}

Since (\ref{u-ni}) satisfies (\ref{PDENeutralInclusions2}), it is left to determine $A_1$ and $\varphi(\rho)$ so that $u$ satisfies the conditions (\ref{conta})-(\ref{transb}) at the interfaces.

Written in ellipsoidal coordinates, the conductivity equation in the coating (\ref{PDENeutralInclusions2}) becomes
\begin{align}
	\label{phicoatingellipcoord}
	\displaystyle 0=\Delta u&=\frac{4g(\rho)}{(\rho-\mu)(\rho-\nu)}\frac{\partial}{\partial\rho}\left[g(\rho)\frac{\partial\Phi}{\partial\rho}\right]\\
&\quad+\frac{4g(\mu)}{(\mu-\rho)(\mu-\nu)}\frac{\partial}{\partial\mu}\left[g(\mu)\frac{\partial\Phi}{\partial\mu}\right]\notag\\
&\quad+\frac{4g(\nu)}{(\nu-\rho)(\nu-\mu)}\frac{\partial}{\partial\nu}\left[g(\nu)\frac{\partial\Phi}{\partial\nu}\right], \notag
\end{align}
where 
\begin{equation}
	\label{g} \displaystyle g(t)=\sqrt{(c_{1}^2+t)(c_{2}^2+t)(c_{3}^2+t)}.
\end{equation}

\begin{remark}
\label{R1}
Observe that $\displaystyle \theta_{1}=\frac{g(\rho_{c})}{g(\rho_{e})}$.
\end{remark}

Using (\ref{u-ni}), (\ref{phicoatingellipcoord}), and the fact that $\Delta x_{1}=0$ we obtain the following second-order differential equation for $\varphi(\rho)$
\begin{equation}
	\label{phicoatingode}
	\displaystyle 0=\frac{d^2\varphi(\rho)}{d\rho^2}+\left[\frac{1}{g(\rho)}\frac{dg(\rho)}{d\rho}+\frac{1}{(c_{1}^2+\rho)}\right]\frac{d\varphi(\rho)}{d\rho}.
\end{equation}
Solving (\ref{phicoatingode}), we obtain
\begin{equation}
	\label{varphi} \displaystyle \varphi(\rho)=A_{2}+B_{2}\int_{\rho_{c}}^{\rho}\frac{1}{(c_{1}^2+\rho)^{\frac{3}{2}}(c_{2}^2+\rho)^{\frac{1}{2}}(c_{3}^2+\rho)^{\frac{1}{2}}}d\rho.
\end{equation}

In what follows, we explain how the unknowns $A_{1}$, $A_{2}$, and $B_{2}$ and $\sigma_{1}^*$ are determined from (\ref{conta}), (\ref{extbound}), (\ref{transa}), and (\ref{transb}).
First, we look at the conditions $u$ must satisfy when $\rho=\rho_c$.  From (\ref{conta}) we have that
\begin{equation}
	\label{interf1} \displaystyle A_{1}=A_{2}+B_{2}\int_{\rho_{c}}^{\rho_{c}}\frac{1}{(c_{1}^2+\rho)^{\frac{3}{2}}(c_{2}^2+\rho)^{\frac{1}{2}}(c_{3}^2+\rho)^{\frac{1}{2}}}d\rho=A_{2},
\end{equation}
and from (\ref{transa}) and (\ref{interf1}), we obtain 
\begin{align}
	\label{interf2} 
	\displaystyle 
B_{2}=\frac{A_{1}g(\rho_{c})(\sigma_{1}\left|A_{1}\right|^{p-2}-\sigma_{2})}{2\sigma_{2}}.
\end{align}

We now look at the conditions that $u$ must satisfy on the outer interface $\rho=\rho_e$.  From (\ref{extbound}) and (\ref{interf1}), we have
\begin{equation}
	\label{extbound1} \displaystyle E=A_{1}+B_{2}\int_{\rho_{c}}^{\rho_{e}}\frac{1}{(c_{1}^2+\rho)^{\frac{3}{2}}(c_{2}^2+\rho)^{\frac{1}{2}}(c_{3}^2+\rho)^{\frac{1}{2}}}d\rho,
\end{equation}
and from (\ref{transb}), we obtain
\begin{align}
	\label{extbound2} \displaystyle 
B_{2}=\frac{Eg(\rho_{e})(\sigma_1^*-\sigma_{2})}{2\sigma_{2}}.
\end{align}

We now introduce the depolarization factors 
\begin{equation}
\label{depfac}
d_{cj}=d_{j}(l_{c1},l_{c2},l_{c3}) \text{, }\hspace{3mm} d_{ej}=d_{j}(l_{e1},l_{e2},l_{e3}) \text{, }\hspace{3mm}j=1,2,3,
\end{equation} 
where 
\begin{equation}
\label{depfacdef}
d_{j}(l_{1},l_{2},l_{3}) = 
\quad\frac{l_{1}l_{2}l_{3}}{2}
\int_{0}^{\infty}\frac{dy}{(l_{j}^2+y)\sqrt{(l_{1}^2+y)(l_{2}^2+y)(l_{3}^2+y)}}
\end{equation} 
is the depolarization factor in direction $j=1,2,3$ of an ellipsoid with semi-axis lenghts $l_{1},l_{2},l_{3}$.  The depolarization factors always sum to unity (see \cite{Milton2002})  
\begin{equation}
\label{dfunit}
d_1+d_2+d_3=1.
\displaystyle 
\end{equation}
Also, observe that $d_{j}(\lambda l_{1},\lambda l_{2},\lambda l_{3})=d_{j}(l_{1},l_{2},l_{3})$ for $\lambda>0$, which means that the depolarization factors are independent of scale.

In terms of these depolarization factors, we have
\begin{align*}
	\displaystyle \int_{\rho_{c}}^{\rho_{e}}\frac{d\rho}{(c_{1}^2+\rho)^{\frac{3}{2}}(c_{2}^2+\rho)^{\frac{1}{2}}(c_{3}^2+\rho)^{\frac{1}{2}}}=\frac{2d_{c1}}{g(\rho_{c})}-\frac{2d_{e1}}{g(\rho_{e})}.
\end{align*}

Rearranging (\ref{extbound2}), we have
\begin{equation}
	\label{extbound4} 
	\displaystyle E=\frac{2B_{2}\sigma_{2}}{g(\rho_{e})(\sigma_1^*-\sigma_{2})}.
\end{equation}
Using (\ref{extbound4}) and (\ref{interf2}) in (\ref{extbound1}), we obtain
\begin{align}
	\label{h2} \displaystyle  \sigma_1^*=\sigma_{2}+\frac{\sigma_{2}\theta_{1}(\sigma_{1}\left|A_{1}\right|^{p-2}-\sigma_{2})}{\sigma_{2}+(\sigma_{1}\left|A_{1}\right|^{p-2}-\sigma_{2})\left[d_{c1}-\theta_{1}d_{e1}\right]}.	
\end{align}
From (\ref{extbound1}), we have
\begin{align}
	\label{extbound5} \displaystyle A_{1}=E-\frac{2B_{2}}{g(\rho_{c})}\left[d_{c1}-\theta_{1}d_{e1}\right]=E-\frac{2B_{2}}{g(\rho_{c})}K,
\end{align}
where $K=d_{c1}-\theta_{1}d_{e1}>0$ is independent of scale.

Using (\ref{extbound5}) in (\ref{interf2}), we obtain the following identity
\begin{align}
	\label{fb2} 
\displaystyle &\sigma_{1}\left|E-\frac{2B_{2}}{g(\rho_{c})}K\right|^{p-2}\left(E-\frac{2B_{2}}{g(\rho_{c})}K\right) \notag\\ &\quad-\sigma_{2}\left(E-\frac{2B_{2}}{g(\rho_{c})}K\right)-\frac{2\sigma_{2}B_{2}}{g(\rho_{c})}=0.
\end{align}

At this point, we consider the function 
\begin{equation}
	\label{fx2} \displaystyle f(x)=\sigma_{1}\left|E-Kx\right|^{p-2}\left(E-Kx\right)-\sigma_{2}\left(E-Kx\right)-\sigma_{2}x.
\end{equation}

Note that we obtain $B_{2}$ if we can prove that $f(x)=0$ has a (unique) solution.  If that is the case, from (\ref{extbound5}) we can obtain $A_{1}$ and from (\ref{h2}) we can get an expression for $\sigma_{1}^*$.

Let us study $f(x)$.  If $E-Kx\geq0$, we have $$f(x)=\sigma_{1}(E-Kx)^{p-1}-\sigma_{2}(E-Kx)-\sigma_2x.$$  Taking the derivative of the $f(x)$, we have $$f'(x)=-K\sigma_{1}(p-1)(E-Kx)^{p-2}+\sigma_{2}(K-1).$$ Note that the first term of $f'(x)$ is negative and the second term is also negative because $K<1$.  To see this, note that by (\ref{dfunit}) and the fact that $K>0$,   
\begin{align*}
K&<
K+(d_{c2}-\theta_{1}d_{e2})+(d_{c3}-\theta_{1}d_{e3})\\
&=(d_{c1}+d_{c2}+d_{c3})-\theta_1(d_{e1}+d_{e2}+d_{e3})\\
&=1-\theta_1=\theta_2<1. 
\end{align*}
Therefore $f(x)$ is a decreasing function.
If $E-Kx<0$, we have $$f(x)=-\sigma_{1}(Kx-E)^{p-1}-\sigma_{2}(E-Kx)-\sigma_2x,$$ and here 
$$f'(x)=-K\sigma_{1}(p-1)(E-Kx)^{p-2}+\sigma_{2}(K-1)$$
is negative for all $x$ so the function $f(x)$ is also decreasing in this case.

Observe that as $x$ approaches $\infty$, the function $f(x)$ approaches $-\infty$ and as $x$ approaches $-\infty$, the function $f(x)$ approaches $\infty$.  Therefore, we conclude that the equation $f(x)=0$ has a unique solution $x_{0}$. 

Moreover, observe that the coefficients of $f(x)$ depend only on $\sigma_{1}$, $\sigma_{2}$, $E$, $K$, and $p$, thus  
\begin{equation}
	\label{ind}
	\displaystyle x_{0}=\frac{2B_{2}}{g(\rho_{c})}=C(\sigma_{1},\sigma_{2},E,K,p).
\end{equation}

Consequently, from (\ref{ind}) and (\ref{extbound5}) we obtain that $\displaystyle A_{1}=E-Kx_{0}$, which together with (\ref{h2}) gives
\begin{equation}
	\label{H*}
	\displaystyle
\sigma_1^* =
\sigma_{2}+\frac{\sigma_{2}\theta_{1}(\sigma_{1}\left|E-\left[d_{c1}-\theta_{1}d_{e1}\right]x_{0}\right|^{p-2}-\sigma_{2})}{\sigma_{2}+(\sigma_{1}\left|E-\left[d_{c1}-\theta_{1}d_{e1}\right]x_{0}\right|^{p-2}-\sigma_{2})\left[d_{c1}-\theta_{1}d_{e1}\right]}.
\end{equation}
Here, we would like to emphasize that (\ref{H*}) shows that $\sigma_{1}^*$ is independent of scale.  In an analogous way, the conductivities in the $x_2$ and $x_3$ directions are obtained and given by similar expresions, also independent of scale.   

\begin{remark}
	If $p=2$, (\ref{fx2}) becomes 
\begin{equation}
	\label{fx2-2} \displaystyle f(x)=\sigma_{1}\left(E-Kx\right)-\sigma_{2}\left(E-Kx\right)-\sigma_{2}x,
\end{equation}
which has a unique root $\displaystyle \bar{x}_{0}=\frac{E(\sigma_1-\sigma_2)}{K(\sigma_1-\sigma_2)+\sigma_2}$. In this case $\sigma_1^*$ (see (\ref{H*})) becomes
\begin{align}
	\label{h2-2-1} \displaystyle  \sigma_1^*=\sigma_{2}+\frac{\sigma_{2}\theta_{1}(\sigma_{1}-\sigma_{2})}{\sigma_{2}+(\sigma_{1}-\sigma_{2})\left[d_{c1}-\theta_{1}d_{e1}\right]}.	
\end{align}
	The conductivities in the $x_2$ and $x_3$ directions are obtain in the same manner and have similar expressions (same results as in Section~7.8 in (\cite{Milton2002})).
\end{remark}

\begin{remark}
	If $c_1=c_2=c_3=c$, we have a sphere.  In this case, (\ref{semiaxis}) becomes
\begin{align}
\label{semiaxis2}
	\displaystyle &\textfrc{l}_{cj}=r_c=\sqrt{c^2+\rho_c} \quad\text{and}\quad \notag\\
	&\textfrc{l}_{ej}=r_e=\sqrt{c^2+\rho_e}, \quad j=1,2,3,
\end{align}
where $r_c$ is the radius of the core of the sphere and $r_e$ the radius of the entire sphere (core and coating).  Here, the volume fraction (\ref{volumefraction}) becomes
\begin{equation}
		\label{volumefraction2}
\theta_{1}=\frac{\textfrc{l}_{c1}\textfrc{l}_{c2}\textfrc{l}_{c3}}{\textfrc{l}_{e1}\textfrc{l}_{e2}\textfrc{l}_{e3}}=\frac{r_c^3}{r_e^3} \text{, and $\theta_2=1-\theta_1$}.
	\end{equation}
The depolarization factors (\ref{depfac}) are all equal and their value is $1/3$, which implies that the integral in (\ref{extbound1}) becomes
\begin{align*}
	\displaystyle \int_{\rho_{c}}^{\rho_{e}}\frac{d\rho}{(c^2+\rho)^{\frac{5}{2}}}=\frac{2}{3}r_c^3\theta_2.
\end{align*}

Therefore we have $\sigma_1^*=\sigma_2^*=\sigma_3^*=\sigma^*$, where
\begin{align}
	\label{H*-sphere}
	\displaystyle
\sigma^*=\sigma_{2}+\frac{3\sigma_{2}\theta_{1}(\sigma_{1}\left|E-\frac{1}{3}\theta_2x_{0}\right|^{p-2}-\sigma_{2})}{3\sigma_{2}+\theta_2(\sigma_{1}\left|E-\frac{1}{3}\theta_2x_{0}\right|^{p-2}-\sigma_{2})},		 
\end{align}
with $x_0$ being the unique and scale-independent solution of
\begin{align}
	\label{fx2-sphere} \displaystyle 
f(x)&=\sigma_{1}\left|E-\frac{1}{3}\theta_2x\right|^{p-2}\left(E-\frac{1}{3}\theta_2x\right) \notag \\
&\quad-\sigma_{2}\left(E-\frac{1}{3}\theta_2x\right)-\sigma_{2}x.
\end{align}
In this way, we recovered the results presented in (\cite{Jimenez2013}).
If $p=2$, we have $$\displaystyle \sigma^*=\sigma_{2}+\frac{3\sigma_{2}\theta_{1}(\sigma_{1}-\sigma_{2})}{3\sigma_{2}+\theta_2(\sigma_{1}-\sigma_{2})},$$ which is the Hashin-Shtrikman formula.
\end{remark}

\bibliographystyle{elsarticle-harv}	

\end{document}